\documentclass[twoside,12pt]{article}
\title{On global defensive $k$-alliances in zero-divisor graph of finite commutative rings}
\date{}
\author{}
\usepackage[all]{xy}
\usepackage{amssymb,amsmath,latexsym}
\usepackage{theorem}
\usepackage{amsfonts}
\usepackage{ulem}
\usepackage{amscd}
\usepackage{amsxtra}
\usepackage{graphicx}
\usepackage{bm}
\usepackage{graphics} 
\usepackage{graphicx} 
\usepackage{pstricks,pst-node} 
\usepackage{tikz} 

\usepackage{tabularx} 
\usepackage{multirow} 
\usepackage{multicol} 
\usepackage{arydshln} 
\usepackage{fancybox} 
\usepackage{multicol} 
\usepackage{array} 
\usepackage{fancybox}

\usepackage{mathrsfs}
\usepackage{array,multirow,makecell}
\setcellgapes{1pt}
\makegapedcells

\newcolumntype{R}[1]{>{\raggedleft\arraybackslash }b{#1}}
\newcolumntype{L}[1]{>{\raggedright\arraybackslash }b{#1}}
\newcolumntype{C}[1]{>{\centering\arraybackslash }b{#1}}

\usetikzlibrary{arrows} 

\newcommand{\field}[1]{\mathbb{#1}}

\newcommand{\Z }{\field{Z}}

\theoremstyle{}\newtheorem{thm}{\bf Theorem}[section]
\theoremstyle{}\newtheorem{cor}[thm]{\bf Corollary}
\theoremstyle{}\newtheorem{lem}[thm]{\bf Lemma}
\theoremstyle{}
\theoremstyle{}\newtheorem{rem}[thm]{\bf Remark}
\theoremstyle{}\newtheorem{pro}[thm]{\bf Proposition}
\theoremstyle{}\newtheorem{exm}[thm]{\bf Example}
\theoremstyle{}\newtheorem{exms}[thm]{\bf Examples}
\theoremstyle{}\newcommand{\cqfd}{\hfill$\square$}

\def\pr{{\parindent0pt {\bf Proof.\ }}}
\def\cqfd
{\hspace{1cm}
\rule{2mm}{2mm}%
\medbreak%
\par%
}

\def\ann{{\rm Ann}}
\def\nil{{\rm Nil}}
\begin{document}
 \thispagestyle{empty}

\maketitle \vspace*{-1.5cm}


\begin{center}
{\large\bf   Driss Bennis, Brahim El Alaoui and Khalid Ouarghi} 
\end{center}
\bigskip 
%
\noindent{\large\bf Abstract.}
The global defensive $k$-alliance is a very well studied notion in graph theory,  it  provides  a method of classification of graphs based on  relations  between members of a particular set of vertices.  In this paper we explore this notion in zero-divisor graph of  commutative rings. The established  results generalize and improve recent work by Muthana and Mamouni who treated a particular case for $k=-1$ known by the global defensive alliance. Various examples  are also provided  which illustrate and  delimit the scope of the established results. \\


\small{\noindent{\bf 2010 Mathematics Subject Classification :}   13M05, 05C25


\section{Introduction} 
Throughout the paper, $R$  will be a commutative ring  with  $1\neq 0$ and $Z(R)$ be  its set of zero-divisors. Let  $x$ be an element  of  $R$, the annihilator of $x$  is defined as  $\ann_R(x):=\{y\in R/ \ xy=0\}$. For an ideal $I$ of $R$, $\sqrt{I}$ means the radical of $I$.  An element $x$ of $R$ is called nilpotent if $x^n=0$ for some positive integers $n$. The set of all nilpotent elements is denoted   $\nil(R):=\sqrt{0}$. A ring $R$ is called reduced if $\nil(R)=\{0\}$. The ring $\Z/n\Z$ of the residues modulo an integer $n$ will be noted by $\Z_n$.  
For a subset $X$ of $R$, we denoted  $X^*=X\setminus \{0\}$.
For any real number $r$, let $\lceil r \rceil$	 (resp.,  $\lfloor r \rfloor$)
denote the ceiling of $r$, that is, the least integer greater than or equal to $r$  (resp., the  floor of $r$, that is  the greatest integer less than or equal to $r$). \\

 Recall that the zero-divisor graph, denoted $\Gamma(R)$,   is a (simple) graph   with vertex set $Z(R)^*,$ and  two distinct vertices $x$ and $y$ are adjacent, if $xy=0$.  We assume the reader has at least a basic familiarity with   zero-divisor graph theory.  For general background   on zero-divisor graph theory, we refer the reader to \cite{AFLL01, ADL, A08, BAD, LEV,  AM07, COY,Axt,DS,SHA}. In this paper, we are interested in studying the global defensive $k$-alliance of zero-divisor graph of finite commutative rings (see Section 2 for the definition of  the  defensive  $k$-alliance of a graph).  It is defined   via    the notions of defensive $k$-alliance and dominating sets (see Section 2).  Several results generalize and improve the recent work of Muthana and  Mamouni  \cite{NA} which focus on the particular where $k=-1$.   

This paper is organized as follows:\\

In Section 2, we recall  the global defensive $k$-alliance  of graphs as well as some notions related to it.
 
In Section 3, we investigate the global defensive $k-$alliance of zero-divisor graph  over local rings. We start by extending \cite[Proposition\ 2.2]{NA} which gives an upper bound on the cardinality of $Z(R)$ in terms of the global defensive alliance number of $\Gamma(R)$. Namely, we will give an upper bound   on the cardinality of $Z(R)$ in terms of the global defensive $k$-alliance number (See Proposition \ref{prop12}). The established upper bound is more optimal than the one given in \cite[Proposition\ 2.2]{NA} as shown by Examples \ref{exam_Z_2timesZ_4andZ_12} and Examples \ref{exam_Z_9andZ_8}. Then, we are interested in studying  the  global defensive $k-$alliance number of the zero-divisor graph over   finite local ring. It seems not an easy task to determine it for any ring.  However, as a second main result,  we  succeed to  compute it for $\Z_{p^{n}}$ for a prime number $p$ and a positive  integer $n$ (see Theorem \ref{thm_Z_p^n}). The  global defensive $k-$alliance number of the  zero-divisor graph  over   finite local ring  with nilpotent maximal ideal  of index $2$ is also given (see Proposition \ref{prop_completegraph}).\\

In  Section 4, we compute the global defensive $k-$alliance number of the zero-divisor graph for some kind of direct products of finite fields. We start by determining  the global defensive $k$-alliance number of the  zero-divisor graph over a direct product of two finite fields (see Theorem \ref{thm2}) and  as a particular case  we fined again \cite[Proposition 2.3]{NA}. Moreover, we get result for global strong defensive alliance which is nothing but the global defensive $0$-alliance (see Section 2).  Determining the global defensive $k$-alliance for the direct product of finite fields  $\prod_{i=1}^{n}F_i$  with $n\geq 3$ a positive integer and $F_i$ a finite field for every $i\in \{1,\ldots,n\}$ is still an open question. However, as main results we determine it for $\Z_2\times \Z_2\times F $ with $|F|\geq2$, and $\Z_2\times F\times K$ with $|K|\geq|F|\geq 3$ (see Theorems \ref{thm_Z_2timesZ_2timesF} and \ref{thm_Z_2timesFtimesK}).\\

Finally,  Section 5 is devoted to study the global defensive $k-$alliance number of zero-divisor graph  for a direct product of  $\Z_2$  and  a finite ring. We start by  giving  an upper and lower bounds for  the global defensive $k$-alliance number of $\Gamma(\Z_2\times R)$ where $R$ is a finite ring (see Theorem \ref{thm4.1}). In \cite[Proposition 2.4]{NA},  Muthana and Mamouni established the equality $\gamma_a(\Gamma(\Z_2\times R))=\lceil\frac{|R|}{2}\rceil$ for a local ring $R$. Here, we   give equalities for     some integers   $k\in [\![1-|R|;1]\!]$ (see Theorems \ref{thm4.2} and \ref{thmfl}).  For a   local ring $R$   with a nilpotent maximal ideal of index $2$, we   improve the inequality of Theorem \ref{thm4.1} and give an equality  for the remaining cases other than the ones studied in Theorems  \ref{thm4.2} and   \ref{thmfl}.


\section{Preliminaries}
 
We  assume some familiarity  with some  basic notions  on    graph theory. Here we deal with the alliance notion of graphs. For reader's convenience, we recall this notion as well as some useful ones including the notion of a dominating set which is a very important notion in graph theory. Ideed, there are plenty of interesting properties related to this notion which still attract the attention of several researchers (see for instance \cite{HHH, TSP, TSP2}). \\

Let $G=(V,E)$ be a finite  simple graph (i.e., a graph without loop or multiple edges). So for two distinct  adjacent  vertices   $x$ and $y$,     we will denote by $x-y$ the edge  between them. 

For a vertex $x\in V,$ the open neighborhood of $x$ is defined as $N(x):=\{y\in V; \ x-y\in E\}$,   and  the closed neighborhood of $x$ is defined by  $N[x]:=N(x)\cup \{x\}$. In general, 
  for a nonempty subset $S\subseteq V$, the open neighborhood of $S$ is defined as $ N(S)=\cup_{x\in S}N(x)$ and its closed neighborhood by  $N[S]=N(S)\cup S$.
  
A set $S$ is dominating if $N[S]=V$.  The minimal cardinality of a dominating set of $G$ is called  the domination number and it is denoted by  $\gamma(G)$.   \\

Next, we give the definition of defensive $k$-alliance (see  \cite{PSS}).  It is worth mentionning that  the  defensive $k$-alliance  has been used to model a variety of applications such as classification problems and distributed protocols (see \cite{Shaf,SX}). To define the defensive $k$-alliance, we first need to recall the definition of degrees of vertices. \\

 The degree of the vertex $x\in V$,  denoted by $deg(x)$,  is the cardinality of its open neighborhood.  Namely,    $deg(x):=|N(x)|$. In general, for every nonempty subset $S\subseteq V$ and every vertex $x\in S$,  we define the degree of $x$ over $S$ as $deg_S(x):=|S\cap N(x)|$. So, $deg_V(x)=deg(x)$. 

 A non-empty set of vertices $S\subseteq V$  is
called a defensive alliance if for every $x\in S$, $|N[x] \cap S|\geq  |N(x) \cap \bar{S}|$,  in  other words, $deg_S(x)+1\geq deg_{\bar{S}}(x)$, where    $\bar{S}=V\setminus S$ (i.e., $\bar{S}$ is the complement of $S$ in $V$). 

A defensive alliance $S$ is called strong if for every vertex $x\in S$, $|N[x]\cap S|> |N(x)\cap \bar{S}|$,  in other words, $deg_S(x)\geq deg_{\bar{S}}(x)$.
In this case we say that every vertex in $S$ is strongly defended. A defensive alliance $S$ is global if it forms a dominating set.\\

The notion defensive alliance   is parametrized in the following sense:\\

First, let us denote by $\Delta(G)$  (resp.,  $\delta(G)$) the maximum  degree (resp., the minimum degree) between all degrees of vertices of $G$. If there is no ambiguity, we simply denote 
 $\Delta$  (resp.,  $\delta$).\\

A non-empty set $S\subseteq V$ is said to be a  defensive $k$-alliance in $G$ with $k$ is an integer in the interval $ [\![-\Delta;\Delta]\!]$,   if for every $x\in S$,  $deg_S(x)\geq deg_{\bar{S}}(x)+k$ or equivalently $deg(x)\geq  2deg_{\bar{S}}(x) +k$ \cite{JIJ}. So, a defensive alliance is nothing but a defensive $(-1)-$alliance and  a strong defensive alliance is nothing but   a defensive $0-$alliance,   as defined in \cite{PSS}. A defensive $0-$alliance which is also known as a cohesive set (see \cite{KR}). 

Notice, that for some graphs, there are some values of $k$ in $ [\![-\Delta;\Delta]\!]$,   such that defensive $k-$alliances do not exist.
For instance, for $k \geq 2$ in the case of the star graph, defensive $k-$alliances do not exist. From the definition of defensive $k-$alliance  we conclude that, in any graph, there are defensive $k-$alliances for $k\in [\![-\Delta;\delta]\!]$.  For instance, a defensive $\delta-$alliance in $G$ is $V$. Moreover, if $x\in V$ is a vertex of minimum degree, $deg(x)=\delta$, then $S=\{x\}$ is a defensive $k-$alliance for every $k\leq -\delta$.\\   
A defensive $k-$alliance set $S$ is called global if it forms a dominating set. The global defensive $k-$alliance number of $G$, denoted by $\gamma_{k}^d(G)$, is the minimum cardinality among all global defensive $k-$alliances in $G$. Clearly, $\gamma_{k+1}^d(G)\geq \gamma_{k}^d(G)\geq \gamma(G)$.
The global defensive $(-1)-$alliance number of $G$ is known as the global alliance number of $G$, denoted by $\gamma_a(G)$,   and the global defensive $0-$alliance number is known as the global strong alliance number,  denoted by $\gamma_{\hat{a}}(G)$ (see \cite{TSM}).


\section{Global defensive $k$-alliances of zero-divisor graph over a local ring}
In this section,  we study the global defensive $k$-alliance of zero-divisor graph over local rings.  We start  by  extending and improving  \cite[Proposition\ 2.2]{NA}.
Before giving the case of local rings, we start by an extension of the first inequality of  \cite[Proposition\ 2.2]{NA}. 
Indeed, the authors  Muthana and Mamouni established the inequality $|Z(R)|\leq \gamma_a(\Gamma(R))^2+\gamma_a(\Gamma(R))+1$.
Here, we extend and ameliorate this inequality to the other global defensive $k-$alliances which provide a more optimal upper bound (see Examples \ref{exam_Z_2timesZ_4andZ_12} and \ref{exam_Z_9andZ_8}).

\begin{lem}\label{lem12}
 For any   finite ring $R$, we have    $ |Z(R)|\leq \min_{k\in [\![-\Delta;\delta]\!]}\{1+\gamma_k^d(\Gamma(R))^2-k\gamma_k^d(\Gamma(R))\}$.\\
Moreover, if  there are   a global defensive $k$-alliance $S=\{x_1,...,x_r\}$ in $\Gamma(R)$ with $r=\gamma_k^d(\Gamma(R))$  and a subset $\Lambda$ of $Z(R)^*$  such that $\Lambda\subset N(x_i)$ for every $x_i\in S$,  then $ |Z(R)|\leq 1+|\Lambda|+r^2-r(k+|\Lambda|)$.
\end{lem}

\pr
Let $k\in [\![-\Delta;\delta]\!]$, and set $r=\gamma_k^d(\Gamma(R))$ and $S=\{x_1,...,x_r\}$ be a global defensive $k$-alliance. We have $N[S]=Z(R)^*$, since $S$ is dominating set. Thus, $\bar{S}\subset \bigcup_{i=1}^r \ann_R(x_i).$ Hence,
\begin{align*}
|\bar{S}|&=|\bar{S}\cap \bigcup_{i=1}^r \ann_R(x_i)|\\&\leq \sum_{i-1}^r|\bar{S}\cap \ann_R(x_i)|\\&\leq \sum_{i=1}^r deg_{\bar{S}}(x_i)\\& \leq \sum_{i=1}^r deg_{S}(x_i)-k \\& \leq \sum_{i=1}^r (r-1)-k=r^2-(1+k)r.
\end{align*}
Then, for every $k\in [\![-\Delta;\delta]\!]$,  $|Z(R)|=1+|S|+|\bar{S}|\leq 1+r^2-kr.$ Hence,  
$|Z(R)|\leq \min_{k\in [\![-\Delta;\delta]\!]}\{1+\gamma_k^d(\Gamma(R))^2-k\gamma_k^d(\Gamma(R))\}.$ Now, if there exists a subset $\Lambda$ of $Z(R)^*$ such that $\Lambda\subset N(x_i)$ for every $x_i\in S.$ Then,
\begin{align*}
 |\bar{S}|&=|\Lambda|+|\bigcup_{i=1}^{r}(\bar{S}\cap \ann_R(x_i))\setminus (\ann_R(x_i)\cap \Lambda)|\\&\leq |\Lambda|+\sum_{i=1}^{r}(|\bar{S}\cap \ann_R(x_i)|-|\ann_R(x_i)\cap \Lambda|)\\&\leq |\Lambda|+\sum_{i=1}^{r}(deg_{\bar{S}}(x_i)-deg_{\Lambda}(x_i))\\&\leq |\Lambda| +\sum_{i=1}^{r}(r-1-k-|\Lambda|)\\&=|\Lambda|+r^2-r-rk-r|\Lambda|
\end{align*}
 and so $|Z(R)|=1+|S|+|\bar{S}|\leq 1+|\Lambda|+r^2-r(k+|\Lambda|)$.
 \cqfd

The following examples show that the bounds in Lemma \ref{lem12}  are optimal. \\
For what follows,  we adopt the following notations: $A_k=1+\gamma_k^d(\Gamma(R))^2-k\gamma_k^d(\Gamma(R))$.

\begin{exms}\label{exam_Z_2timesZ_4andZ_12}
 	\begin{enumerate}
 	\item Let $R=\Z_{12}$. The zero-divisor graph over  this ring is illustrated in Figure \ref{fig1}. And different values of  the global defensive $k$-alliance number and $A_k$ are presented in Table $\ref{Tab1}$.
 \begin{figure}[ht]
 	\centering
 	\includegraphics[scale=0.5]{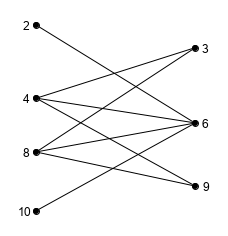}
 	\caption{$\Gamma(R)$}
 	\label{fig1}
 \end{figure}
 \begin{table} 
 	\begin{center}
 		\begin{tabular}{|C{0.5cm}|C{3cm}|C{4cm}|}
 			\hline $k$ & $\gamma_k^d(\Gamma(R))$ &  $A_k$ \\
 			\hline  $-4$ & $2$ & $13$  \\
 			\hline  $-3$ & $2$ & $11$  \\
 			\hline  $-2$ & $2$ & $9$   \\
 			\hline  $-1$ & $3$ & $13$  \\
 			\hline  $0$ & $4$ & $17$   \\
 			\hline  $1$ & $5$ & $21$   \\
 			\hline 
 		\end{tabular} 
 		\label{Tab1}
 	\end{center}
 	\caption{Values of   $\gamma_k^d(\Gamma(R))$ and $A_k$} 
 \end{table}
 This example gives an   upper bound of the cardinality of $Z(R)$ which is smaller than the upper bound giving in \cite[Proposition 2.2]{NA}. Namely,  $|Z(R)|< \min_{k\in [\![-4;1]\!]}\{ A_k\}=9<\gamma_a(\Gamma(R))^2+\gamma_a(\Gamma(R))+1=13$.
\item Consider the ring $R=\Z_2\times \Z_4$. The zero-divisor graph over  this ring is illustrated in Figure \ref{fig2}. And different values of  the global defensive $k$-alliance number and $A_k$ are presented in Table $\ref{Tab1.2}$.  
 		\begin{figure}[ht]
 			\centering
 			\includegraphics[scale=0.5]{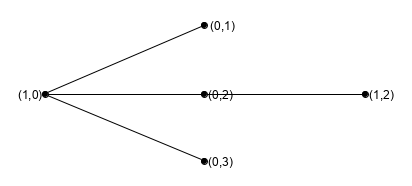}   
 			\caption{$\Gamma(R)$}   
 			\label{fig2}
 		\end{figure}

 \begin{table}
 		\begin{center}
 			\begin{tabular}{|C{0.5cm}|C{3cm}|C{4cm}|}
 				\hline $k$ & $\gamma_k^d(\Gamma(R))$ &  $A_k$ \\
 				\hline  $-3$ & $2$ & $11$  \\
 				\hline  $-2$ & $2$ & $9$   \\
 				\hline  $-1$ & $2$ & $7$  \\
 				\hline  $0$ & $3$ & $10$   \\
 				\hline  $1$ & $4$ & $7$   \\
 				\hline 
 			\end{tabular}
 			\label{Tab1.2}
 		\end{center}
 	\caption{Values of $\gamma_k^d(\Gamma(R))$ and $A_k$  }
 \end{table}

This ring provides an example satisfying the second inequality for $k=-1$. Namely,   we have  the global defensive $(-1)$-alliance $\{(1,0),(1,2)\}$ and the set $\Lambda =\{(0,2)\}$ satisfying the condition of Lemma \ref{lem12}, that is,  $\Lambda \subseteq N((1,0))\cap N((1,2))$.  Then, $|Z(R)|\leq 1+|\Lambda| +\gamma_{-1}^d(\Gamma(R))^2-\gamma_{-1}^d(\Gamma(R))(-1+|\Lambda|)=6 $.\\
 Notice that  $|Z(R)|< \min_{k\in [\![-3;1]\!]}\{ A_k\}=\gamma_a(\Gamma(R))^2+\gamma_a(\Gamma(R))+1=7$. 
\end{enumerate}
\end{exms}

Now, we give our first desired result.
 
\begin{pro}\label{prop12}
Let $R$ be a  finite local ring. Then, $|Z(R)|\leq \max\{\\ \min_{k\in [\![-\Delta;\delta]\!]} \{2\gamma_k^d(\Gamma(R))-k\}, \ \min_{k\in [\![-\Delta;\delta]\!]}\{2+\gamma_k^d(\Gamma(R))^2-(k+1)\gamma_k^d(\Gamma(R))\}\}$.
\end{pro}

\pr
Let $k\in [\![-\Delta;\delta]\!]$ and set $r=\gamma_k^d(\Gamma(R))$ and $S=\{x_1,...,x_r\}$ be a global defensive $k$-alliance. We have $N[S]=Z(R)^*$, since $S$ is dominating set. Thus, $\bar{S}\subset \bigcup_{i=1}^r \ann_R(x_i)$ and since   $R$ is a finite local ring,  there exists $x\in Z(R)^*$ such that  its  maximal ideal $M$ is  $Z(R)=\ann_R(x)$. 
If $x\notin S$, then   by taking $\Lambda=\{x\}$ and using Lemma \ref{lem12},  we get the following inequality  for every $k\in [\![-\Delta;\delta]\!]$:
  $$|Z(R)|=1+|S|+|\bar{S}|\leq 2+r^2-(k+1)r.$$ 
If $x\in S$,  then  $deg_S(x)\geq deg_{\bar{S}}(x)+k.$ That is $|S|-1\geq |\bar{S}|+k$. Then, for every $k\in [\![-\Delta;\delta]\!]$, $|Z(R)|=1+|S|+|\bar{S}|\leq 2r-k$.
Hence $|Z(R)|\leq max\{\min_{k\in [\![-\Delta;\delta]\!]}\{2r-k\}, \min_{k\in [\![-\Delta;\delta]\!]}\{2+r^2-(k+1)r\}\}$.
\cqfd

We give  examples proving that the bounds given  in Proposition \ref{prop12}  are sharp.\\ 
For what follows, we adopt the notations: $B_k=2\gamma_k^d(\Gamma(R))-k$ and $C_k=2+\gamma_k^d(\Gamma(R))^2-(k+1)\gamma_k^d(\Gamma(R))$.

\begin{exms}\label{exam_Z_9andZ_8}
\begin{enumerate}
\item Let  $R=\Z_9$. The  different values of  the global defensive $k$-alliance number and $A_k$ are presented in Table  $\ref{Tab1.3}$. 
\begin{table}
\begin{center}
\begin{tabular}{|C{0.5cm}|C{2cm}|C{2.5cm}|C{2.5cm}|}
\hline $k$ & $\gamma_k^d(\Gamma(R))$  & $B_k$ & $C_k$ \\
\hline  $-1$ & $1$  & $3$ & $3$ \\
\hline  $0$ & $2$  & $4$ & $4$  \\
\hline  $1$ & $2$  & $3$ & $2$  \\
\hline 
\end{tabular}
\label{Tab1.3}
\end{center}
\caption{Values of  $\gamma_k^d(\Gamma(R))$, $B_k$ and  $C_k$   }
\end{table}
Then, 
\begin{equation*}
 \left\{
\begin{array}{ll}
& \min_{k\in [\![-1;1]\!]}\{B_k\}=3, \\&
\min_{k\in [\![-1;1]\!]}\{C_k\}=2,
\end{array}
\right.
\end{equation*} 
and so   
$|Z(R)|=\max\{\min_{k\in \{-1,0,1\}}\{B_k\}, \ \min_{k\in \{-1,0,1\}}\{C_k\}\}=3$.

\item Consider the ring $R=\Z_8$. The zero-divisor graph over  this ring is illustrated in Figure \ref{fig3}. And different values of  the global defensive $k$-alliance number and $A_k$ are presented in Table $\ref{Tab1.4}$.   
\begin{figure}[ht]
\centering
\includegraphics[scale=0.5]{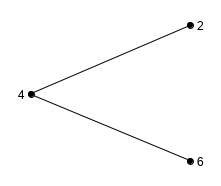}   
\caption{$\Gamma(R)$}   
\label{fig3}
\end{figure}
\begin{table}
\begin{center}
\begin{tabular}{|C{0.5cm}|C{2cm}|C{2.5cm}|C{2.5cm}|C{2.5cm}|}
\hline $k$ & $\gamma_k^d(\Gamma(R))$   & $B_k$ & $C_k$ \\
\hline  $-2$ & $1$  & $4$ & $4$  \\
\hline  $-1$ & $2$  & $5$ & $6$ \\
\hline  $0$ & $2$  & $4$ & $4$  \\
\hline  $1$ & $3$ & $5$ & $5$  \\
\hline 
\end{tabular}
\label{Tab1.4}
\end{center}
\caption{Values of $\gamma_k^d(\Gamma(R))$,  $B_k$ and  $C_k$ }
\end{table}
Then, 
\begin{equation*}
\left\{
\begin{array}{ll}
& \min_{k\in [\![-2;1]\!]}\{B_k\}=4, \\&
\min_{k\in [\![-2;1]\!]}\{C_k\}=4,
\end{array}
\right.
\end{equation*}

and so  $|Z(R)|=\max\{\min_{[\![-2;1]\!]}\{B_k\}, \ \min_{[\![-2;1]\!]}\{C_k\}\}=4$.
\end{enumerate}
\end{exms}

These examples  showing also  that $\min_{[\![-\Delta;\delta]\!]} B_k$ and $\min_{[\![-\Delta;\delta]\!]} C_k$ are not comparable in sens that $\min_{[\![-\Delta;\delta]\!]} C_k<\min_{[\![-\Delta;\delta]\!]} B_k$ for certain rings but they are not for others.

\begin{rem}
 Notice that for   a finite local ring $R$, $\gamma_k^d(\Gamma(R))\geq 2+k$ if and only if  $2\gamma_{k}^d(\Gamma(R))-k\leq  2+(\gamma_{k}^d(\Gamma(R)))^2-(k+1)\gamma_k^d(\Gamma(R))$. So, For every $k\in [\![-\Delta;\delta]\!],$
  	\begin{equation*}
  	|Z(R)|\leq   \left\{
  	\begin{array}{ll}
  	& 2\gamma_{k}^d(\Gamma(R))-k \text{\ \ \ \ \ \ \ \ \ \ \ \ \ \ \ \ \ \ \ \ \ \ \  \ \ \ \ \ \ \  if } \gamma_k^d(\Gamma(R))\leq 2+k,\\&
  2+(\gamma_{k}^d(\Gamma(R)))^2-(k+1)\gamma_k^d(\Gamma(R)) \text{\  \ \ \ \ } otherwise.
  	\end{array}
  	\right.
  	\end{equation*}    
\end{rem}

It is not clear how to determine global defensive $k-$alliances of any finite local ring. But we can determine it for $\Z_{p^{n}}$ for a prime number $p$ and an integer $n$. 
However, we know that for a finite local ring $(R,M)$, $\Gamma(R)$ is complete if and only if $Z(R)=M$ with $M^2=0$,   \cite[Theorem 2.8]{ADL}. And we know from \cite{JJ} that the global defensive $k$-alliance number of a complete graph is determined as follows: for every $k\in [\![1-n;n-1]\!]$, $\gamma_k^d(K_n)=\lceil\frac{n+k+1}{2}\rceil$. So, we have the following consequence for this simple case.

\begin{pro}\label{prop_completegraph}
Let $R$ be a finite local ring such that  its maximal ideal $M$ is nilpotent of index $2$. Then,  for every $k\in[\![2-|M|;|M|-2]\!]$,  $\gamma_k^d(\Gamma(R))=\lceil\frac{|M|+k}{2}\rceil$. 
\end{pro}

Idealization  can be used to give a family of examples of rings whose zero-divisor graph is complete. Recall  that  the idealization of an $R-$module $M$ called also the trivial extension of $R$ by $M$, denoted by $R(+)M$, is the commutative ring $R\times M$ with the following addition and multiplication: $(a,n)+(b,m)=(a+b,n+m)$ and $(a,n)(b,m)=(ab,am+bn)$  for every   $(a,n),(b,m)\in  R(+)M$.

\begin{exm}
   Let  $n$ be a positive integer and $p$ be a prime number. Then, $\Z_p(+)(\Z_{p})^n$ is a local finite ring of maximal ideal  $0(+)(\Z_{p})^n.$  We have $(0(+)(\Z_{p})^n)^2=0$ and  so 
   $\Gamma(\Z_p(+)(\Z_{p})^n)$ is a complete graph. Then,  $\gamma_{k}^d(\Z_p(+)(\Z_{p})^n)=\lceil \frac{p^n+k}{2}\rceil$ for every $k\in [\![1-p^n;p^n-1]\!].$
\end{exm}

Now, we give the main result of this section. We determine the global defensive $k$-alliance number of the zero-divisor graph $\Gamma(\Z_{p^n})$.  If $n=2$, then the maximal ideal $Z(\Z_{p^2})=<p>$  is  nilpotent of index $2$ and so by  Proposition \ref{prop_completegraph} we have $\gamma_{k}^d(\Gamma(\Z_{p^2}))=\lceil\frac{p+k}{2}\rceil$ for every $k\in [\![2-p;p-2]\!]$. For $n\geq 3$, we have the following theorem.

\begin{thm}\label{thm_Z_p^n}
	Let $p$ be a prime number and $n\geq 2$ be an integer. Then, for every $k\in [\![2-p^{n-1};p-1]\!]$,   $\gamma_{k}^d(\Gamma(\Z_{p^n}))=\lceil\frac{p^{n-1}+k}{2}\rceil.$ 
\end{thm}
\pr 
We have $Z(\Z_{p^n})=\{\overline{ap}/ \ 0\leq a< p^{n-1} \}$ and $|Z(\Z_{p^n})|=p^{n-1}$. For each $1\leq r\leq n-1$, we define $A_r=\{\overline{ap^r}/ \ p\ \text{does not divide}\ a\}$ and so   $Z(\Z_{p^n})=\bigcup_{r=1}^{r=n-1}A_r$. There are two cases to discuss:\\
\textbf{Case} $p=2:$ Let $k\in [\![2-2^{n-1};1]\!]$ and set $S=S_1\cup \{2^{n-1}\}$ such that $S_1\subseteq A_1$ with $|S_1|=\lceil \frac{2^{n-1}+k}{2}\rceil -1$. 	 We have $deg_S(\bar{2})=1$ and $deg_{\bar{S}}(\bar{2})+k=0+k$. Then, $deg_{\bar{S}}(\bar{2})+k\leq deg_S(\bar{2})$. And, $deg_S(2^{n-1})=|S_1|=\lceil \frac{2^{n-1}+k}{2}\rceil -1$ and $deg_{\bar{S}}(2^{n-1})+k=|\bar{S}|+k=|Z(\Z_{p^n})|-1-|S|+k=2^{n-1}+k-1- \lceil\frac{2^{n-1}+k}{2}\rceil\leq \lceil\frac{2^{n-1}+k}{2}\rceil-1$. So,  $deg_{\bar{S}}(2^{n-1})+k \leq deg_S(2^{n-1})$. Hence, $ S $ is a global defensive $k$-alliance of cardinality $\lceil\frac{2^{n-1}+k}{2}\rceil$.\\
Now, let $S$ be a global defensive $k$-alliance of minimal cardinality $\gamma_{k}^d(\Gamma(\Z_{2^n}))$. 
If $2^{n-1}\notin S$. Then, $A_1\subseteq S$ and so $deg_S(\bar{2})\geq deg_{\bar{S}}(\bar{2})+k=1+k$. In cases $k=0$ and  $k=1$,  we get a contradiction.  Hence, $2^{n-1}\in S$ and so  $deg_S(2^{n-1})\geq deg_{\bar{S}}(2^{n-1})+k$.  Then,  $|S|-1\geq |\bar{S}|+k=|Z(\Z_{2^n})|-1-|S|+k$ and so  $|S|\geq \frac{2^{n-1}+k}{2}$. Thus,   $|S|\geq \lceil\frac{2^{n-1}+k}{2}\rceil$.  Hence,  $\gamma_{k}^d(\Gamma(\Z_{2^n}))= \lceil\frac{2^{n-1}+k}{2}\rceil.$\\	
\textbf{Case} $p\geq 3:$  Let $k\in [\![2-p^{n-1};p-1]\!]$ and $S$ be a global defensive $k$-alliance of minimal cardinality $\gamma_{k}^d(\Gamma(\Z_{p^n})).$  If $A_{n-1}\cap S=\emptyset.$ Then, $A_1\subseteq S$ and so $deg_S(\bar{p})\geq deg_{\bar{S}}(\bar{p})+k=|A_{n-1}|+k=p-1+k$ a contradiction when $k\in [\![2-p;p-1]\!].$ Then, $A_{n-1}\cap S\neq\emptyset$ then,  there exists $z\in A_{n-1}\cap S$ such that $deg_S(z)\geq deg_{\bar{S}}(z)+k$ and so $|S|-1\geq |Z(\Z_{p^n})^*|-|S|+k.$ Thus, $|S|\geq \frac{p^{n-1}+k}{2}.$ Hence, $\gamma_{k}^d(\Gamma(\Z_{p^n}))\geq \lceil \frac{p^{n-1}+k}{2}\rceil.$ \\ 
 Now, to prove that the cardinality of $S$ is less than or equal  to $\lceil \frac{p^{n-1}+k}{2}\rceil$ we have to  prove that there exists a global defensive $k$-alliance of cardinality $\lceil \frac{p^{n-1}+k}{2}\rceil$ for every $k\in [\![2-p^{n-1};p-1]\!].$ So we have the following sub-cases:\\
\underline{\textbf{Sub-case}  $k\in [\![2-p^{n-1};2(p-1)-p^{n-1}]\!]$}:\\
Let   $S\subseteq A_{n-1}$ such that $|S|=\lceil \frac{p^{n-1}+k}{2}\rceil.$ It is clear that $S$ is a dominating set.  Let $x\in S$, then   $deg_S(x)=|S|-1$  and  $deg_{\bar{S}}(x)+k=|Z(\Z_{p^n})|-1-|S|+k=p^{n-1}+k-1-\lceil \frac{p^{n-1}+k}{2}\rceil \leq \lceil \frac{p^{n-1}+k}{2}\rceil -1$ and so  $deg_S(x)\geq deg_{\bar{S}}(x)+k.$ Hence, $S$ is a global defensive $k$-alliance of cardinality $|S|=\lceil \frac{p^{n-1}+k}{2}\rceil$.\\
\underline{\textbf{Sub-case}  $k\in [\![2(p-1)-p^{n-1}+1;2(p^2-1)-p^{n-1}]\!]$}:\\ 
Let $S\subseteq A_{n-1}\cup A_{n-2}$ with $S\cap A_{n-1}\neq \emptyset$ and $S\cap A_{n-2}\neq \emptyset$  such that $|S|=\lceil \frac{p^{n-1}+k}{2}\rceil.$ Clearly, $S$ is a dominating set. Let $x\in S\cap A_{n-1},$ we have $deg_S(x)=|S|-1$ and $deg_{\bar{S}}(x)+k=|Z(p^{n-1})|-1-|S|+k=p^{n-1}+k-\lceil \frac{p^{n-1}+k}{2}\rceil +k\leq |S|-1$ and so $deg_S(x)\geq deg_{\bar{S}}(x)+k$. Let $x\in S\cap A_{n-2},$  
we have $deg_S(x)=|S|-1$ and $deg_{\bar{S}}(x)+k=|Z(\Z_{p^n})|-1-(|S|+|A_1|)+k=p^{n-1}-1-|S|-|A_1| +k=p^{n-1}+k-\lceil \frac{p^{n-1}+k}{2}\rceil -|A_1| -1\leq \lceil \frac{p^{n-1}+k}{2}\rceil-1$ and so $deg_S(x)\geq deg_{\bar{S}}(x)+k.$  Hence, $S$ is a global defensive $k$-alliance of cardinality $|S|=\lceil \frac{p^{n-1}+k}{2}\rceil$.\\
 \underline{\textbf{Sub-case}  $k\in [\![2(p^2-1)-p^{n-1}+1;2(p^3-1)-p^{n-1}]\!]$}:\\Let $S\subseteq A_{n-1}\cup A_{n-2}\cup A_{n-3}$ with $S\cap A_{n-1}\neq \emptyset$, $S\cap A_{n-2}\neq \emptyset$ and $S\cap A_{n-3}\neq \emptyset$  such that $|S|=\lceil \frac{p^{n-1}+k}{2}\rceil.$ Clearly, $S$ is a dominating set. Let $x\in S\cap A_{n-1}$,   we have $deg_S(x)=|S|-1$ and $deg_{\bar{S}}(x)+k=|Z(\Z_{p^{n}})|-1-|S|+k=p^{n-1}+k-\lceil \frac{p^{n-1}+k}{2}\rceil +k\leq |S|-1$ and so $deg_S(x)\geq deg_{\bar{S}}(x)+k$. Let $x\in S\cap A_{n-2}$  we have $deg_S(x)=|S|-1$ and $deg_{\bar{S}}(x)+k=|Z(\Z_{p^n})|-1-(|S|+|A_1|)+k=p^{n-1}-1-|S|-|A_1| +k=p^{n-1}+k-\lceil \frac{p^{n-1}+k}{2}\rceil -|A_1| -1\leq \lceil \frac{p^{n-1}+k}{2}\rceil-1$ and so $deg_S(x)\geq deg_{\bar{S}}(x)+k$. Now, let $x\in S\cap A_{n-3}$,  then  $deg_S(x)=|S|-1$ and $deg_{\bar{S}}(x)+k=|Z(\Z_{p^n})|-1-(|S|+|A_1|+|A_2|)+k=p^{n-1}-1-|S|-|A_1|-|A_2| +k=p^{n-1}+k-\lceil \frac{p^{n-1}+k}{2}\rceil -|A_1|-|A_2| -1\leq \lceil \frac{p^{n-1}+k}{2}\rceil-1$ and so $deg_S(x)\geq deg_{\bar{S}}(x)+k.$  Hence, $S$ is a global defensive $k$-alliance of cardinality $|S|=\lceil \frac{p^{n-1}+k}{2}\rceil$.\\
 \underline{\textbf{Sub-cases}  $k\in [\![2(p^{\alpha-1}-1)-p^{n-1}+1;2(p^{\alpha}-1)-p^{n-1}]\!]$ for  $  3 \leq  \alpha\leq \frac{n}{2}$}:\\
 Let $S\subseteq A_{n-1}\cup A_{n-2}\cup \dots \cup A_{n-\alpha}$ with $S\cap A_{n-1} \neq \emptyset$,   $S\cap A_{n-2} \neq \emptyset$, \dots,  $S\cap A_{n-\alpha} \neq \emptyset.$ such that $|S|= \lceil \frac{p^{n-1}+k}{2}\rceil.$ Thus, similarly to the  previous sub-case, we     prove that $S$ is a global defensive $k$-alliance of cardinality $|S|=\lceil \frac{p^{n-1}+k}{2}\rceil$.\\ 
\underline{\textbf{Sub-case}  $k\in  [\![2(p^{\frac{n}{2}}-1)-p^{n-1}+1;2(p^{\frac{n}{2}+1}-1)-p^{n-1}]\!]$}:\\
Let $S\subseteq A_{n-1}\cup A_{n-2}\cup \dots \cup A_{\frac{n}{2}-1}$ with $S\cap A_{n-1} \neq \emptyset$,   $S\cap A_{n-2} \neq \emptyset$, \dots,  $S\cap A_{\frac{n}{2}-1} \neq \emptyset$ and 
\begin{equation}\label{equ_1}
		 \left\{
		 \begin{array}{ll}
		 & |S\cap A_{\frac{n}{2}-1}|\leq \lfloor\frac{|A_1|+|A_2|+\dots+|A_{\frac{n}{2}-1}|}{2}\rfloor,\\&
		 |S\cap A_{\frac{n}{2}-1}|+|S\cap A_{\frac{n}{2}}|\leq \lfloor\frac{|A_1|+|A_2|+\dots+|A_{\frac{n}{2}}|}{2}\rfloor
		 \end{array}
		 \right.
\end{equation}
such that $|S|= \lceil \frac{p^{n-1}+k}{2}\rceil$. It  is clear that $S$ is  dominating set. Let $x\in S\cap A_\beta$ with $\frac{n}{2}+1\leq \beta \leq n-1$, we have  $deg_S(x)=|S|-1$ and $deg_{\bar{S}}(x)+k=|Z(\Z_{p^n})|-1-|S|-(|A_1|+\dots+|A_{n-\beta-1}|)$. Then,  $deg_S(x)\geq deg_{\bar{S}}(x)+k$.
 Let $x\in S\cap A_{\frac{n}{2}}$ we have $deg_S(x)=|S|-1-|S\cap A_{\frac{n}{2}-1}|$ and $deg_{\bar{S}}(x)+k=|Z(\Z_{p^n})|-1-|S|-(|A_1|+\dots +|A_{\frac{n}{2}-2}|+|\bar{S}\cap A_{\frac{n}{2}-1}|)+k\leq |S|-1-(|A_1|+\dots +|A_{\frac{n}{2}-2}|+|\bar{S}\cap A_{\frac{n}{2}-1}|)$. Then, by  (\ref{equ_1}),   $deg_S(x)\geq deg_{\bar{S}}(x)+k$. 
 Let $x\in S\cap A_{\frac{n}{2}}$,  we have $deg_S(x)=|S|-1-(|S\cap A_{\frac{n}{2}}|+|S\cap A_{\frac{n}{2}-1}|)$ and $deg_{\bar{S}}(x)+k=|Z(\Z_{p^{n}})|-1-|S|-(|A_1|+|A_2|+\dots + |A_{\frac{n}{2}-2}|+|\bar{S}\cap A_{\frac{n}{2}-1}|+|\bar{S}\cap A_{\frac{n}{2}}|)\leq |S|-1-(|A_1|+|A_2|+\dots + |A_{\frac{n}{2}-2}|+|\bar{S}\cap A_{\frac{n}{2}-1}|+|\bar{S}\cap A_{\frac{n}{2}}|) $. 
 So, by   (\ref{equ_1}),  $deg_S(x)\geq deg_{\bar{S}}(x)+k$. Hence, $S$ is a global defensive $k-$alliance of cardinality  $|S|=\lceil \frac{p^{n-1}+k}{2}\rceil$.\\
 \underline{\textbf{Sub-cases}  $k\in  [\![2(p^{\frac{n}{2}+\alpha-1}-1)-p^{n-1}+1;2(p^{\frac{n}{2}+\alpha}-1)-p^{n-1}]\!]$ for $2\leq \alpha\leq \frac{n}{2}-2$}:\\
Let $S\subseteq A_{n-1}\cup A_{n-2}\cup \dots \cup A_{\frac{n}{2}-\alpha}$ with $S\cap A_{n-1}\neq \emptyset$, $S\cap A_{n-2}\neq \emptyset$, \dots,  $S\cap A_{\frac{n}{2}-\alpha}\neq \emptyset$ and  

\begin{equation*}
\left\{
\begin{array}{ll}
& |S\cap A_{\frac{n}{2}-\alpha}|\leq \lfloor\frac{|A_1|+\dots +|A_{\frac{n}{2}-\alpha}|}{2}\rfloor,\\&
|S\cap A_{\frac{n}{2}-\alpha}|+|S\cap A_{\frac{n}{2}-\alpha+1}|\leq\lfloor\frac{|A_1|+\dots +|A_{\frac{n}{2}-\alpha+1}|}{2}\rfloor,\\&
....\ \ \ \ \ \ .... \ \ \ \ \  ....... \ \ \ \ \ \ ....\\&

|S\cap A_{\frac{n}{2}-\alpha}|+|S\cap A_{\frac{n}{2}-\alpha+1}|+\dots + |S\cap A_{\frac{n}{2}}|\leq \lfloor\frac{|A_1|+|A_2|+\dots+|A_{\frac{n}{2}}|}{2}\rfloor
\end{array}
\right.
\end{equation*}
such that $|S|=\lceil\frac{p^{n-1}+k}{2}\rceil$. Then, similarly to  the previous  sub-case, we     prove that   $S$ is a global defensive $k-$alliance.\\
\underline{\textbf{Sub-case}  $k\in [\![2(p^{n-2}-1)-p^{n-1}+1;p-1]\!]$}:\\
Let $S\subseteq A_{n-1}\cup A_{n-2}\cup \dots \cup A_1$ with $S\cap A_{n-1} \neq \emptyset$,   $S\cap A_{n-2} \neq \emptyset$, \dots,  $S\cap A_{1} \neq \emptyset$ and 
\begin{equation*}
	 \left\{
	 \begin{array}{ll}
	 & |S\cap A_1|\leq \lfloor\frac{|A_1|}{2}\rfloor,\\&
	 |S\cap A_1|+|S\cap A_2|\leq \lfloor\frac{|A_1|+|A_2|}{2}\rfloor,\\&
	 ....\ \ \ \ \ \ .... \ \ \ \ \  ....... \ \ \ \ \ \ ....\\&
	 
	 |S\cap A_1|+|S\cap A_2|+\dots + |S\cap A_{\frac{n}{2}}|\leq \lfloor\frac{|A_1|+|A_2|+\dots+|A_{\frac{n}{2}}|}{2}\rfloor
	 \end{array}
	 \right.
\end{equation*}
such that $|S|= \lceil \frac{p^{n-1}+k}{2}\rceil$. So, it is not difficult to verify that $S$ is a global defensive $k$-alliance of cardinality $\lceil \frac{p^{n-1}+k}{2}\rceil$. \\
 Hence, for every $k\in [\![2-p^{n-1};p-1]\!]$,  $\gamma_{k}^d(\Gamma(\Z_{p^n}))= \lceil \frac{p^{n-1}+k}{2}\rceil$.
\cqfd


\section{Global defensive $k$-alliances of zero-divisor graph of some kind of direct product of finite fields}
In this section, we study  the global defensive $k-$alliance number of zero-divisor graphs over some kind of direct products of finite fields.\\

The global defensive alliance  and the global strong alliance  numbers for a star graph and   complete bipartite graph are determined as follow:

\begin{pro} [\cite{TSM},  Proposition 3]
The following assertions hold: 
	\begin{enumerate}
		\item $\gamma_a(K_{1,s})=\lfloor \frac{s}{2}\rfloor +1$.
		\item $\gamma_a(K_{r,s})=\lfloor\frac{r}{2}\rfloor + \lfloor \frac{s}{2}\rfloor$, if $r,s\geq 2.$
		\item $\gamma_{\hat{a}}(K_{r,s})=\lceil\frac{r}{2}\rceil + \lceil\frac{s}{2}\rceil.$
	\end{enumerate}
\end{pro}

In the following result,  we determine  the global defensive $k$-alliance number,   $k\in [\![-\Delta;\delta]\!]$, of zero-divisor graphs over a  direct product of two finite fields.

\begin{thm}\label{thm2}
	Let $F$ and $K$ be two  finite fields. Then,
	\begin{enumerate}
		\item  If $|F|=2$ $(F\cong \Z_2)$  and $|K|\geq 2.$ Then, for every $k\in [\![1-|K|;1]\!]$,   $$\gamma_k^d(\Gamma(F\times K))=\lceil \frac{|K|+k+1}{2} \rceil$$ 
		\item If $|K|\geq|F|\geq 3,$
		 \begin{equation*}
		\gamma_k^d(\Gamma(F\times K))=  \left\{
		\begin{array}{ll}
		& 2 \text{\ \ \ \ \ \ \ \ \ \ \ \ \ \  \ \ \ \ \ \ \  if } k={1-|K|},\\&
		\lceil\frac{|K|+k+1}{2} \rceil \ \text{\  \ \ \ \ \ \ \ \ \ if } k\in [\![2-|K|;3-|F|]\!], \\&
		\lfloor\frac{|F|+k}{2}\rfloor + \lfloor\frac{|K|+k}{2}\rfloor  \text{ if } k\in [\![4-|F|;|F|-1]\!].
		\end{array}
		\right.
		\end{equation*} 
	\end{enumerate}
\end{thm}

\pr
1.  Let $k\in [\![1-|K|;1]\!] $ and let $S$ be a global defensive $k-$alliance of cardinality $r=\gamma_k^d(\Gamma(F\times K))$. 
If $(1,0)\notin S.$ Then, $(0,u)\in S$ for all $u\in K^*.$ Thus,  for $k=0\  (resp., \  k=1)$, we have $deg_S((0,u))=0\geq deg_{\bar{S}}((0,u))=1\ (resp.,\ deg_S((0,u))=0\geq deg_{\bar{S}}((0,u))+1=2),$ a contradiction. Consequently, $(1,0)\in S$ and so  $deg_S((1,0))\geq deg_{\bar{S}}((1,0))+k $. That is  $|S\cap \{0\}\times K^*|\geq |\bar{S}\cap \{0\}\times K^*|+k.$ Then,  $|S\cap \{0\}\times K^*|\geq |\{0\}\times K^*|-|\bar{S}\cap \{0\}\times K^*|+k.$ Thus, $|S\cap \{0\}\times K^*|\geq \frac{|K|-1+k}{2}.$  Hence,  $ r \geq |S\cap \{0\}\times K^*|+1\geq \frac{|K|+k+1}{2}$ and so  $ r \geq \lceil\frac{|K|+k+1}{2}\rceil.$
Now, set $S=\{(1,0)\}\cup \{0\}\times S_1,$ such that $S_1\subseteq K^*$ and $|S_1|=\lceil\frac{|K|+k+1}{2}\rceil-1.$ Clearly, $S$ is dominating set. 
We have $deg_S((1,0))=|S_1|=\lceil \frac{|K|+k+1}{2}\rceil -1$  and
$deg_{\bar{S}}((1,0))+k=|\{0\}\times K^*|-|S_1|+k =|K|+k-\lceil\frac{|K|+k+1}{2}\rceil$. Then,   $deg_{\bar{S}}((1,0))+k\leq deg_S((1,0))$. Let $u\in K^*$,  then 
$deg_S((0,u))=1\geq k= deg_{\bar{S}}(0,u)+k.$
Hence, $S$ is a global defensive $k$-alliance and so $r\leq |S|=\lceil\frac{|K|+k+1}{2}\rceil.$ Then, by the first part $r=\lceil\frac{|K|+k+1}{2}\rceil.$\\
2. {\textbf{Case} $k=1-|K|$}:\\ Since  every dominating set in $\Gamma(F\times K )$  has a cardinality greater than or equal $2.$ Then, every global defensive $k$-alliance $S$ has cardinality $|S|\geq 2.$\\
Now, set $S=\{(0,v),(u,0)\}$ such that $v\in K^*$ and $u\in F^*.$ It is clear that $S$ is a dominating set. \\
We have $deg_S((0,v))=1$ and $deg_{\bar{S}}((0,v))+1-|K|=|F|-1+1-|K|=|F|-|K|.$ Then, $deg_S((0,v))\geq deg_{\bar{S}}((0,v))+k.$ Moreover, $deg_S((u,0))=1$ and $deg_{\bar{S}}((u,0))+1-|K|=|K|-2+1-|K|=-1.$ Thus, $deg_S((u,0))\geq deg_{\bar{S}}((u,0))+k$. So,  $S$ is a global defensive $k$-alliance of cardinality $|S|=2$. Hence, $r=\gamma_{k}^d(\Gamma(F\times K))=|S|=2$.\\   
{\textbf{Case} $k\in [\![2-|K|;3-|F|]\!]$}:\\
Set $S=\{(u,0)\}\cup \{0\}\times S_1$ such that $u\in F^*$,  $S_1\subseteq K^*$ and $|S_1|=\lceil\frac{|K|+k+1}{2}\rceil-1$. Clearly, $S$ is dominating set. 
We have $deg_S((u,0))=|S_1|=\lceil \frac{|K|+k+1}{2}\rceil -1$  and
$deg_{\bar{S}}((u,0))+k=|\{0\}\times K^*|-|S_1|+k = |K|+k-\lceil\frac{|K|+k+1}{2}\rceil\leq \lceil \frac{|K|+k+1}{2}\rceil-1$. Then,   $deg_{\bar{S}}((u,0))+k\leq deg_S((u,0))$.\\
Let $v\in S_1.$ We have, $deg_S((0,v))=1$ and $ deg_{\bar{S}}(0,v)+k=|\bar{S}\cap F^*\times \{0\}|+k=|F|-2+k$. Thus,  $deg_S((0,v))\geq deg_{\bar{S}}(0,v)+k$. Hence, $S$ is a global defensive $k$-alliance of cardinality $|S|=\lceil\frac{|K|+k+1}{2}\rceil$.\\ 
Now, let $S$ be a global defensive $k-$alliance of minimal cardinality $r=\gamma_k^d(\Gamma(\Z_2\times F))$. Since $S$ is a dominating set, it must be contain element of the form $ (0,v)$ and $(u,0)$ with $v\in K^*$ and $u\in F^*$. 
Suppose that there exist $(u_1,0),(u_2,0)\in S$ with $u_1\neq u_2$. Then, for every $(0,v)\in S$  we have $deg_S((0,v))=|S\cap F^*\times \{0\}|\geq 2$. On the other hand $deg_S((u_1,0))\geq deg_{\bar{S}}((u_1,0))+k$ and so $|S\cap \{0\}\times K^*|\geq \frac{|K|+k-1}{2}$.  Hence, $ \lceil\frac{|K|+k+1}{2}\rceil\geq |S|\geq \frac{|K|+k-1}{2}+2=\frac{|K|+k+1}{2}+1$ a contradiction.  Then, $u_1=u_2$ and so $deg_S((0,v))=1$ for every $(0,v)\in S$. Hence,  $ \lceil\frac{|K|+k+1}{2}\rceil\geq |S|\geq \frac{|K|+k-1}{2}+1$ and so $r=|S|=\lceil\frac{|K|+k+1}{2}\rceil$.\\
{\textbf{Case} $k\in [\![4-|F|;|F|-1]\!]$}:\\ Let $S=S_1\times \{0\}\cup \{0\}\times S_2$ such that $S_1\subseteq F^*$ and $S_2\subseteq K^*$ with $|S_1|=\lfloor \frac{|F|+k}{2}\rfloor$ and  $|S_2|=\lfloor \frac{|K|+k}{2}\rfloor.$ Clearly,  $S$ is a dominating set.  We have,  $deg_{S}((u,0))=|S_2|=\lfloor \frac{|K|+k}{2}\rfloor$ and $deg_{\bar{S}}((u,0))+k=|\bar{S}\cap \{0\}\times K^*|+k=|\{0\}\times K^*|-|S_2|+k=|K|+k-1-\lfloor\frac{|K|+k}{2}\rfloor\leq   \frac{|K|+k}{2}-\frac{1}{2}\leq \lfloor\frac{|K|+k}{2}\rfloor.$ Then, $deg_S((u,0))\leq deg_{\bar{S}}((u,0))+k$. 
Let $v\in S_2$, then  $deg_S((0,v))=|S_1|=\lfloor \frac{|F|+k}{2}\rfloor$  and $deg_{S}((0,v))+k=|\bar{S}\cap F^*\times \{0\}|+k=|F^*\times \{0\}|-|S\cap F^*\times \{0\}|+k=|F|+k-1-\lfloor \frac{|F|+k}{2}\rfloor.$ Thus, $deg_{\bar{S}}((0,v))+k\leq  \frac{|F|+k}{2} -\frac{1}{2}\leq deg_S((0,v))$. 
Hence, $S$ is global defensive $k$-alliance.\\
Now, Let $S$ be global defensive $k$-alliance of cardinality $r=\gamma_k^d(\Gamma(F\times K))$. 
Assume that for all $u\in F^*$, $(u,0)\notin S.$ Then, $(0,v)\in S$  for all $v\in K^*$ since $S$ is a dominating set and so $deg_S((0,v))\geq deg_{\bar{S}}((0,v))+k.$ Then, $0\geq |\bar{S}\cap F^*\times \{0\}|+k=|F|-1+k\geq 1$  a contradiction. We get also a contradiction when we suppose that $(0,v)$ is not in $S$ for all  $v$ in $K^*$.  Hence, $(0,v), (u,0)\in S$ for some $u\in F^*$ and $v\in K^*.$ Then, $deg_S((u,0))\geq deg_{\bar{S}}((u,0))+k$, that is  $|S\cap \{0\}\times K^*|\geq |\bar{S}\cap \{0\}\times K^*|+k.$ Thus,  $|S\cap \{0\}\times K^*|\geq | \{0\}\times K^*|-|S\cap \{0\}\times K^*|+k $  and so $|S\cap \{0\}\times K^*|\geq \frac{|K|-1+k}{2}.$ Similarly, we have  $deg_S((0,v))\geq deg_{\bar{S}}((0,v))+k$ and so  $	|S\cap F^*\times \{0\}|\geq \frac{|F|-1+k}{2}.$ 	Hence, $|S|=|S\cap \{0\}\times K^*| + |S\cap F^*\times \{0\}|\geq\lfloor\frac{|F|+k}{2}\rfloor+\lfloor\frac{|K|+k}{2}\rfloor$ and so  $r=\lfloor\frac{|F|+k}{2}\rfloor+\lfloor\frac{|K|+k}{2}\rfloor$.
\cqfd
As a particular case of  Theorem \ref{thm2}, we fined again \cite[Proposition 2.3]{NA}. Moreover, we get result of global strong defensive alliance which corresponds to  global defensive $0$-alliance.

\begin{cor}
Let $F$ and $K$ be finite fields. Then,
\begin{enumerate}
\item $\gamma_a(\Gamma(\Z_2\times F))=\lfloor\frac{|F|-1}{2}\rfloor + 1$,  \cite[Proposition 2.3]{NA}.
\item  $\gamma_{\hat{a}}(\Gamma(\Z_2\times F))=\lceil\frac{|F|}{2}\rceil + 1$.
\item $\gamma_a(\Gamma(F\times K))=\lfloor\frac{|F|-1}{2}\rfloor + \lfloor\frac{|K|-1}{2}\rfloor$, if $|F|, |K|\geq 3,$ \cite[Proposition 2.3]{NA}.
\item $\gamma_{\hat{a}}(\Gamma(F\times K))=\lceil\frac{|F|-1}{2}\rceil + \lceil\frac{|K|-1}{2}\rceil$.
\end{enumerate}
\end{cor}

\begin{cor}
Let $p$ and $q$ be two  prime numbers. Then, 
\begin{enumerate}
	\item \textbf{Case} $p=2$ and $q\geq 2:$  For every $k\in [\![1-q;1]\!],$  $\gamma_{k}^d(\Gamma(\Z_p\times \Z_q))=\lceil\frac{q+k+1}{2}\rceil.$
	\item \textbf{Case}  $3\leq p\leq q:$  
	 \begin{equation*}
	\gamma_k^d(\Gamma(\Z_p\times \Z_q))=  \left\{
	\begin{array}{ll}
	& 2 \text{\hspace{2.25cm}  if } k={1-q},\\&
	\lceil\frac{q+k+1}{2} \rceil \ \text{\hspace{0.9cm} if } k\in [\![2-q;3-p]\!], \\&
	\lfloor\frac{p+k}{2}\rfloor + \lfloor\frac{q+k}{2}\rfloor  \text{ if } k\in [\![4-p;p-1]\!].
	\end{array}
	\right.
	\end{equation*}
\end{enumerate}
\end{cor}

In \cite[Theorem 2.6]{NA},  Muthana and Mamouni established the equality $\gamma_a(\Gamma(\Z_2\times \Z_2\times F  ))=|F|$ for a finite field  $F$ of cardinality $|F|\geq 3$.  Here, we give equalities for $|F|\geq 2$ and for all  $k\in [\![1-2|F|;1]\!]$.

\begin{thm}\label{thm_Z_2timesZ_2timesF}
Let $F$ be a finite field. Then,
	\begin{equation*}
	\gamma_k^d(\Gamma(\Z_2\times \Z_2\times F))=  \left\{
	\begin{array}{ll}
	& 	3 \text{\  \ \ \ \ \ \ \ \ \ \ \ \  \ if } k\in [\![1-2|F|;3-2|F|]\!],\\&
	 |F|+\lceil\frac{1+k}{2}\rceil \text{  if } k\in [\![4-2|F|;1]\!].
	\end{array}
	\right.
	\end{equation*}
\end{thm}

\pr 
\textbf{Case} $k\in [\![1-2|F|;3-2|F|]\!]$:  There is no dominating set in $\Gamma(\Z_2\times \Z_2\times F)$  of  cardinality smaller  than or equal to two. On the other hand, we consider the set $S=\{(0,0,z),(1,0,0),(0,1,0)\}$ with $z\in F^*$ and we can  prove easily that $S$ is a global defensive $k$-alliance with $k\in [\![1-2|F|;3-2|F|]\!]$.\\
\textbf{Case} $k\in [\![4-2|F|;1]\!]$: Let $S=\{0\}\times \{0\}\times S_1\cup \{(1,0,0),(0,1,0)\}$ with $S_1$ is a subset of $F^*$ and $|S_1|=|F|+\lceil\frac{1+k}{2}\rceil -2$. Clearly, S is a dominating set. We have,  $deg_S((0,1,0))=1+|S_1|=|F|+\lceil\frac{1+k}{2}\rceil -1$ and $deg_{\bar{S}}((0,1,0))+k=|\Z_2\times \{0\}\times F|-1-(|\{0\}\times \{0\}\times S_1|+1)+k=|F|-\lceil\frac{1+k}{2}\rceil+k\leq |F|+\lceil\frac{1+k}{2}\rceil-1$.  Then, $deg_S((0,1,0))\geq deg_{\bar{S}}((0,1,0))+k$. Similarly, $deg_S((1,0,0))\geq deg_{\bar{S}}((1,0,0))+k$. 
Let $z\in S_1$, then  $deg_S((0,0,z))=2$ and $deg_{\bar{S}}((0,0,z))+k=|\Z_2\times \Z_2 \times \{0\}|-1-2+k=1+k$.  Thus, $deg_S((0,0,z))\geq  deg_{\bar{S}}((0,0,z))+k$. 
Hence, $S$ is a global defensive $k$-alliance of cardinality $|S|=|F|+\lceil\frac{1+k}{2}\rceil$.\\
Now, let $S$ be a global defensive $k$-alliance of minimal cardinality $r=\gamma_k^d(\Gamma(\Z_2\times\Z_2\times F))$. 
From the first part, we have $r\leq  |F|+\lceil\frac{1+k}{2}\rceil$. 
Suppose that $(1,0,0)\notin S$, then  $(0,1,z)\in S$ for all $z\in F^*$,  thus   $deg_S((0,1,z))\geq deg_{\bar{S}}((0,1,z))+k$ and so $|S\cap \Z_2 \times\{0\}\times \{0\}|=0\geq 1+k$ which is not true  for $k=0,1$. Hence, $(1,0,0)\in S$ and so $deg_S((1,0,0))\geq deg_{\bar{S}}((1,0,0)) +k$, then $|S\cap \{0\}\times \Z_2\times F|\geq |F|+\frac{k-1}{2}$. Thus, $|F|+\frac{k+1}{2}\leq |S|\leq |F|+\lceil\frac{1+k}{2}\rceil$. Hence, $r=|F|+\lceil\frac{1+k}{2}\rceil$. 		
\cqfd

\begin{cor} 
Let $p$ be a prime number. Then, 
\begin{equation*}
     \gamma_k^d(\Gamma(\Z_2\times \Z_2\times \Z_p))=  \left\{
     \begin{array}{ll}
     & 	3 \text{\hspace{1.5cm} if } k\in [\![1-2p;3-2p]\!],\\&
     p+\lceil\frac{1+k}{2}\rceil \text{  if } k\in [\![4-2p;1]\!].
     \end{array}
     \right.
\end{equation*}
\end{cor}

In the following theorem we extend the equality  $\gamma_a(\Gamma(\Z_2\times F \times K))=\lceil\frac{|F||K|}{2}\rceil$ giving in \cite[Theorem 2.8]{NA} to the other cases $k\in [\![1-|K||F|;1]\!]$.

\begin{thm}\label{thm_Z_2timesFtimesK}
Let $F$ and $K$  be two  finite fields with $|K|\geq |F|\geq 3$.  Then,
	\begin{equation*}
	\gamma_k^d(\Gamma(\Z_2\times F \times K))=  \left\{
	\begin{array}{ll}
	& 3 \text{ \hspace{2.7cm} if }  k\in [\![1-|K||F|;5-|K||F|]\!],\\&
	\lceil\frac{|F||K|+k+1}{2}\rceil \text{\hspace{1cm} if } k\in [\![6-|F||K|;1]\!].
	\end{array}
	\right.
	\end{equation*}
\end{thm}

\pr
\textbf{Case} $k\in [\![1-|K||F|;5-|K||F|]\!]$:
Let $S$ be a global defensive $k$-alliance of minimal cardinality $r=\gamma_{k}^d(\Gamma(\Z_2\times F\times K))$.  Since $S$ is a dominating set,  then the cardinality of $S$ must be greater than or equal $3$. On the other hand,  set $S=\{(1,0,0),(0,u,0),(0,0,v)\}$ with $u\in F^*$ and $v\in K^*.$ It is clear that $S$ is a dominating set.  We have,  $deg_S((1,0,0))=2$ and $deg_{\bar{S}}((1,0,0))+k=(|K^*||F^*|+|K^*|-1+|F^*|-1)+k=|F||K|-3+k\leq |F||K|-3+5-|F||K|=2$. Then, $deg_{\bar{S}}((1,0,0))+k\leq deg_S((1,0,0))$. 
We have,  $deg_S((0,u,0))=2$ and $deg_{\bar{S}}((0,u,0))+k=|K^*|+|K^*|-1+k=2|K|-3+k\leq 2|K|-3+5-|K||F|\leq 2+|K|(2-|F|)$. Then, $deg_{\bar{S}}((0,u,0))+k\leq deg_S((0,u,0))$ (since $|F|>2$). 
We have,  $deg_S((0,0,v))=2$ and $deg_{\bar{S}}((0,0,v))+k=|F^*|+|F^*|-1+k=2|F|-3+k\leq 2|F|-3+5-|K||F|\leq 2+|F|(2-|K|)$. Then, $deg_{\bar{S}}((0,0,v))+k\leq deg_S((0,0,v))$.  Thus, $S$ is global defensive $k$-alliance of cardinality $|S|=3$.  Hence, $r=3$. \\
\textbf{Case}  $k\in [\![6-|K||F|;5-2|K|]\!]$: Let $S=\{(1,0,0), (0,u,0), (0,0,v)\}\cup S_1$ with $S_1\subseteq \{0\}\times F^*\times K^*$ such that $|S_1|=\lceil\frac{|F||K|+k+1}{2}\rceil-3$. Clearly,  $S$ is a dominating set.
We have,  $deg_S((1,0,0))=|S_1|+2=\lceil\frac{|K||F|+k+1}{2}\rceil -1$ and $deg_{\bar{S}}((1,0,0))+k=|F||K|-3-|S_1|+k=|F||K|-\lceil\frac{|K||F|+k+1}{2}\rceil+k\leq \lceil\frac{|K||F|+k+1}{2}\rceil -1$. So, $deg_S((1,0,0))\geq deg_{\bar{s}}((1,0,0)) +k$. We have,  $deg_S((0,u,0))=2$ and $deg_{\bar{S}}((0,u,0))+k=2|K|-3+k$ and so $deg_S((0,u,0))\geq deg_{\bar{S}}((0,u,0))+k$. Also, $deg_S((0,0,v))=2$ and $deg_{\bar{S}}((0,0,v))+k=2|F|-3+k$. Then,  $deg_S((0,0,v))\geq deg_{\bar{S}}((0,0,v))+k$. Let $(0,v_1,v_2)\in S_1$, we have $deg_S((0,v_1,v_2))=1\geq k=deg_{\bar{S}}((0,v_1,v_2))+k$. Hence, $S$ is a global defensive $k-$alliance of cardinality $|S|=\lceil\frac{|F||K|+k+1}{2}\rceil$.\\ Now, let $S$ be a global defensive $k-$alliance of minimal cardinality $\gamma_{k}^d(\Gamma(\Z_2\times F\times K))$. If $(1,0,0)\notin S$, then $\{0\}\times F^*\times K^*\subseteq S$ and also $S$  contains at least two vertices from $\{0\}\times F^*\times \{0\}\cup \{0\}\times \{0\}\times K^*\cup \{1\}\times \{0\}\times K^*\cup \{1\}\times F^*\times \{0\}$. Then, $|S|\geq |F^*||K^*|+2$ or from the first part,  $|S|\leq \lceil\frac{|F||K|+k+1}{2}\rceil$ and so $ |F^*||K^*|+2\leq  \lceil\frac{|F||K|+k+1}{2}\rceil$ which is  not true for every $k\in [\![6-|K||F|;5-2|K|]\!]$. Hence,  $(1,0,0)\in S$  and so $deg_S((1,0,0))\geq deg_{\bar{S}}((1,0,0))+k$, then $|S\cap \{0\}\times F\times K|\geq |\{0\}\times F\times K|-|S\cap \{0\}\times F\times K|+k-1$. Thus, $|S\cap \{0\}\times F\times K|\geq \frac{|F||K|+k-1}{2}$ and so $|S|\geq \frac{|F||K|+k-1}{2}+1=\frac{|F||K|+k+1}{2}$. 
\\
\textbf{Case}   $k\in [\![6-2|K|;5-2|F|]\!]$: Let $S=\{(1,0,0),(0,u,0)\}\cup \{0\}\times \{0\}\times S_1\cup S_2$ with $S_1\subseteq K^*$ and $S_2\subseteq \{0\}\times F^*\times K^*$ such that $|S_1|=\lceil\frac{2|K|+k+1}{2}\rceil-2$ and $|S_2|=\lceil\frac{|K||F|+k+1}{2}\rceil-\lceil\frac{2|K|+k+1}{2}\rceil$. Clearly, $S$ is a dominating set. We have,  $deg_S((1,0,0))=|S_1|+|S_2|+1=\lceil\frac{|K||F|+k+1}{2}\rceil-1$ and $deg_{\bar{S}}((1,0,0))+k=|\{0\}\times F^*\times K^*|-|S_2|+|K^*|-|S_1|+|F^*|-1+k=|F||K|+k-\lceil\frac{|K||F|+k+1}{2}\rceil\leq \lceil\frac{|K||F|+k+1}{2}\rceil-1$. So,  $deg_S((1,0,0))\geq  deg_{\bar{S}}((1,0,0))+k$. We have, $deg_S((0,u,0))=|S_1|+1=\lceil\frac{2|K|+k+1}{2}\rceil -1$ and $deg_{\bar{S}}((0,u,0))+k=|K^*|+|K^*|-|S_1|+k=2|K|+k-\lceil\frac{2|K|+k+1}{2}\rceil\leq \lceil\frac{2|K|+k+1}{2}\rceil-1$. Then,  $deg_{\bar{S}}((0,u,0))+k\leq  deg_S((0,u,0))$. Let $v\in S_1$, we have $deg_S((0,0,v))=2$ and $deg_{\bar{S}}((0,0,v))+k=|F^*|-1+|F^*|+k$. So, $deg_{\bar{S}}((0,0,v))+k\leq deg_S((0,0,v))$. Let $(0,v_1,v_2)\in S_2$, then $deg_S((0,v_1,v_2))=1$ and $deg_{\bar{S}}((0,v_1,v_2))+k=0+k=k$. So, $deg_{\bar{S}}((0,v_1,v_2))+k\leq  deg_S((0,v_1,v_2))$. Hence, $S$ is a global defensive $k-$alliance of cardinality $|S|=\lceil\frac{|F||K|+k+1}{2}\rceil$.\\ Now, let $S$ be a global defensive $k-$alliance of minimal cardinality $\gamma_{k}^d(\Gamma(\Z_2\times F\times K))$. If $(1,0,0)\notin S$, then $\{0\}\times F^*\times K^*\subset S$ and also  $S$  contains at least two vertices from $\{0\}\times F^*\times \{0\}\cup \{0\}\times \{0\}\times K^*\cup \{1\}\times \{0\}\times K^*\cup \{1\}\times F^*\times \{0\}$. Then, $|S|\geq |F^*||K^*|+2$ and   $|S|\leq \lceil\frac{|F||K|+k+1}{2}\rceil$ from the first part. So,  $ |F^*||K^*|+2\leq  \lceil\frac{|F||K|+k+1}{2}\rceil$ which is  not true since $k\in [\![6-2|K|;5-2|F|]\!]$. Hence,  $(1,0,0)\in S$  and so $deg_S((1,0,0))\geq deg_{\bar{S}}((1,0,0))+k$, then $|S\cap \{0\}\times F\times K|\geq |\{0\}\times F\times K|-|S\cap \{0\}\times F\times K|+k-1$. Thus,  $|S\cap \{0\}\times F\times K|\geq \frac{|F||K|+k-1}{2}$ and so $|S|\geq \frac{|F||K|+k-1}{2}+1=\frac{|F||K|+k+1}{2}$. Hence, $|S|=\lceil\frac{|F||K|+k+1}{2}\rceil$.\\ 
\textbf{Case} $k\in [\![6-2|F|;1]\!]$:
Set $S=\{(1,0,0)\}\cup \{0\}\times S_1\times \{0\}\cup \{0\}\times \{0\} \times S_2\cup S_3$ with  $S_1\subseteq F^*$, $S_2\subseteq K^*$, and $S_3\subseteq \{0\}\times F^*\times K^*$ such that $|S_1|=\lceil\frac{2|F|+k+1}{2}\rceil-2 $ $|S_2|=\lceil\frac{2|K|+k+1}{2}\rceil-2$ and $|S_3|=\lceil\frac{|F||K|+k+1}{2}\rceil-\lceil\frac{2|F|+k+1}{2}\rceil- \lceil\frac{2|K|+k+1}{2}\rceil+3$.  Clearly, $S$ is a dominating set. \\
We have $deg_S((1,0,0))=|S_1|+|S_2|+|S_3|= \lceil\frac{|F||K|+k+1}{2}\rceil-1$ and $deg_{\bar{S}}((1,0,0))+k= |F^*||K^*|-|S_3|+|F^*|-|S_1|+|K^*|-|S_2|+k=(|F^*||K^*|+|F^*|+|K^*|)-(|S_3|+|S_2|+|S_1|)+k=|F||K|+k-\lceil\frac{|F||K|+k+1}{2}\rceil\leq \lceil\frac{|F||K|+k+1}{2}\rceil-1$. Then, $deg_S((1,0,0))\geq deg_{\bar{S}}((1,0,0))+k$. 
Let $u\in S_1,$ then $deg_S((0,u,0))=1+|S_2|=\lceil\frac{2|K|+k+1}{2}\rceil-1$ and $deg_{\bar{S}}((0,u,0))+k=|K^*|-|S_2|+|K^*|+k=2|K|+k-\lceil\frac{2|K|+k+1}{2}\rceil\leq \lceil\frac{2|K|+k+1}{2}\rceil-1 $. So,  $deg_S((0,u,0))\geq deg_{\bar{S}}((0,u,0))+k$.
Let $v\in S_2,$ then $deg_S((0,0,v))=1+|S_1|=\lceil\frac{2|F|+k+1}{2}\rceil-1$ and $deg_{\bar{S}}((0,0,v))+k=|F^*|-|S_1|+|F^*|+k=2|F|+k-\lceil\frac{2|F|+k+1}{2}\rceil\leq \lceil\frac{2|F|+k+1}{2}\rceil-1$ So,  $deg_S((0,0,v))\geq deg_{\bar{S}}((0,0,v))+k$.
Let $(0,u,v)\in S_3,$ we have $deg_S((0,u,v))=1$ and $deg_{\bar{S}}((0,u,v))+k=0+k$. Thus, $deg_S((0,u,v))\geq deg_{\bar{S}}((0,u,v))+k$.
Hence, $S$ is a global defensive $k$-alliance of cardinality $|S|=\lceil\frac{|F||K|+k+1}{2}\rceil$.\\
Now, let $S$ be a global defensive $k$-alliance of minimal cardinality $r=\gamma_{k}^d(\Gamma(\Z_2\times F\times K))$. 
If $(1,0,0)\notin S,$ then $\{0\}\times F^*\times K^*\subset S$ and so $deg_S((0,u,v))=0\geq deg_{\bar{S}}((0,u,v))+k=1+k$ which is not true for $k\in \{0,1\}$. Hence, $(1,0,0)\in S$ and so $deg_S((1,0,0))\geq deg_{\bar{S}}((1,0,0))+k$, then $|S\cap \{0\}\times F\times K|\geq |\{0\}\times F\times K|-1-|S\cap \{0\}\times F\times K|+k$. Thus,  $|S\cap \{0\}\times F\times K|\geq \frac{|F||K|+k-1}{2}$ and so $|S|\geq \frac{|F||K|+k-1}{2}+1=\frac{|F||K|+k+1}{2}$. Since,  $|S|\leq \lceil\frac{|F||K|+k+1}{2}\rceil$ (by the first part), $r=|S|= \lceil\frac{|F||K|+k+1}{2}\rceil$.	
\cqfd

\begin{cor}
	Let $p\geq 3$ and $q\geq 3$ be two prime numbers. Then,  
\begin{equation*}
	\gamma_k^d(\Gamma(\Z_2\times \Z_p \times \Z_q))=  \left\{
	\begin{array}{ll}
	& 3 \text{ \hspace{2.2cm} if }  k\in [\![1-pq;5-pq]\!],\\&
	\lceil\frac{pq+k+1}{2}\rceil \text{\hspace{1cm} if } k\in [\![6-pq;1]\!].
	\end{array}
	\right.
\end{equation*}
\end{cor}

\section{ Global defensive $k$-alliances of zero-divisor graph of the  direct  product of $\Z_2$ and a finite  ring}

The purpose  of this section is to study the global defensive $k-$alliances of zero-divisor graphs over the  direct product of  $\Z_2$  and a  finite ring.

We start by  giving  an upper and lower bounds of  the global defensive $k$-alliance number of $\Gamma(\Z_2\times R)$ for every $k\in [\![1-|R|;1]\!]$  where $R$ is a finite ring.
\begin{thm}\label{thm4.1}
	Let $R$ be a finite ring. Then, for every $k\in [\![1-|R|;1]\!]$,
	  $$\lceil\frac{|R|+k+1}{2}\rceil\leq \gamma_k^d(\Gamma(\Z_2\times R))\leq\lceil\frac{|R|+2|Z(R)|+k-1}{2}\rceil.$$
\end{thm}
\pr
If $R$ is a finite field. Then, $\gamma_k^d(\Gamma(\Z_2\times R))= \lceil\frac{|R|+k+1}{2}\rceil$, by Theorem \ref{thm2}. Thus, we can assume that $R$ is not a finite  field. Let $S$ be a global defensive $k$-alliance of minimal cardinality $r=\gamma_{k}^d(\Gamma(\Z_2\times R))$. 
If $(1,0)\notin S$,   then $\{0\}\times U(R)\subseteq S.$ So, for $y\in U(R)$ we have $deg_S((0,y))\geq deg_{\bar{S}}((0,y))+k$ and so $k\leq -1$ which is not true for $k=0$ and $k=1$. Hence,    $(1,0)$ must be inside of $S$ and so  $deg_S((1,0))\geq deg_{\bar{S}}((1,0))+k.$ Then, $|S\cup \{0\}\times R^*|\geq |\bar{S}\cup \{0\}\times R^*|+k$ and so $|S\cup \{0\}\times R^*|\geq \frac{|R|+k-1}{2}.$ Then, $|S|\geq \frac{|R|+k+1}{2}$. Hence, $r=|S|\geq \lceil\frac{|R|+k+1}{2}\rceil.$ It is easy to check the other bound ( you have just to see that  we have to add $|Z(R)^*|$  element  from $\Z_2\times Z(R)^*$  to have the property of domination), that is $\gamma_k^d(\Gamma(\Z_2\times R))\leq\lceil\frac{|R|+k+1}{2}\rceil+|Z(R)|-1=\lceil\frac{|R|+2|Z(R)|+k-1}{2}\rceil$. 
\cqfd

In \cite[Proposition 2.4]{NA},  Muthana and Mamouni established the equality $\gamma_a(\Gamma(\Z_2\times R))=\lceil\frac{|R|}{2}\rceil$ for a local ring $R$. Here, we give equalities for some other integers   $k\in [\![1-|R|;1]\!]$. The following result present the cases $k=0$ and $k=1$.

\begin{thm}\label{thm4.2}
 Let $R$ be a finite local ring which is not a field. Then,
  \begin{equation*}
 \gamma_{k}^d(\Gamma(\Z_2\times R))=  \left\{
 \begin{array}{ll}
 & \lceil \frac{|R|+1}{2} \rceil \text{ \hspace{1cm} if } k=0, \\&
 \lceil\frac{|R|}{2}\rceil+2 \text{\hspace{1cm}if } k=1.
 \end{array}
 \right.
 \end{equation*}
\end{thm}
\pr 
\textbf{Case} $k=0$:  We have,  $|Z(R)|\leq \lceil\frac{|R|}{2}\rceil$.  Then, we can consider a set $S_1\subseteq R^*$  such that $Z(R)^*\subseteq S_1$ and $|S_1|=\lceil \frac{|R|}{2}\rceil$. 
Set  $S=\{(1,0)\}\cup \{0\}\times S_1$ and let  prove that is a global strong alliance.
Clearly $S$ is dominating set. We have,  $deg_S((1,0))=|\{0\}\times S_1|=|S_1|=\lceil\frac{|R|}{2}\rceil$  and $deg_{\bar{S}}((1,0))=|\bar{S}\cap \{0\}\times R^*|=|\{0\}\times R^*|-|S_1|=|R|-\lceil\frac{|R|}{2}\rceil-1\leq \lceil \frac{|R|}{2}\rceil -1$. Then,   $deg_S((1,0))\geq deg_{\bar{S}}((1,0))$.
Let $z\in Z(R)^*$. If $z^2=0$, then
$deg_S((0,z))=|S\cap\Z_2\times \ann_R(z)|=|\{(1,0)\}\cup \{0\}\times \ann_R(z)|-2=|\ann_R(z)|-1$ and $deg_{\bar{S}}((0,z))=|\bar{S}\cap \Z_2\times \ann_R(z)|=|\{1\}\times \ann_R(z)|-1=|\ann_R(z)|-1$. Thus,  $deg_S((0,z)\geq deg_{\bar{S}}((0,z))$.  If $z^2\neq 0$. Then, 
$deg_S((0,z))=|S\cap\Z_2\times \ann_R(z)|=|\{(1,0)\}\cup \{0\}\times \ann_R(z)|-1=|\ann_R(z)|$  and $deg_{\bar{S}}((0,z))=|\bar{S}\cap \Z_2\times \ann_R(z)|=|\{1\}\times \ann_R(z)|-1=|\ann_R(z)|-1$. Thus,  $deg_S((0,z)\geq deg_{\bar{S}}((0,z))$.  Let $u\in S_1\setminus Z(R)$. Then $deg_S((0,u))=1\geq 0=deg_{\bar{S}}((0,u))$. Hence,  $S$ is a  global strong alliance set  and so  $\gamma_0^d(\Gamma(\Z_2\times R))\leq \lceil \frac{|R|}{2} \rceil +1$.\\
Now, let $S$ be a global strong alliance of minimal cardinality  $r=\gamma_0^d(\Z_2\times R)$. 
Suppose that $(1,0)\notin S$. Then,  for every element $u\in U(R)$, $(0,u)\in S,$ since $S$ is a dominating set. Or,  $S$ is a strong alliance set, then $deg_S((0,u))=0\geq deg_{\bar{S}}((0,u))=1$ a contradiction. Hence, $(1,0)\in S$ and so  $deg_S((1,0))\geq deg_{\bar{S}}((1,0))$ which implies that  $|S\cap \{0\}\times R^*|\geq |\bar{S}\cap \{0\}\times R^*|=|\{0\}\times R^*|-|S\cap \{0\}\times R^*|$ and so $|S\cap \{0\}\times R^*|\geq \frac{|R|-1}{2}.$ Thus, $|S|\geq \frac{|R|}{2}-\frac{1}{2}+1$. Hence,  $\frac{|R|}{2}+\frac{1}{2}\leq r\leq \lceil\frac{|R|}{2}\rceil +1 $  which implies that $\gamma_0^d(\Gamma(\Z_2\times R))= \lceil \frac{|R|+1}{2} \rceil$.\\
\textbf{Case} $k=1$: Consider the  set $S_1\subseteq R^*$ such that $Z(R)^*\subseteq S_1$ and $|S_1|=\lceil \frac{|R|}{2}\rceil$.  Let $x\in Z(R)^*$ such that $\ann_R(x)=Z(R)$ ($R$ is local ring) and set $S=\{(1,0),(1,x)\}\cup \{0\}\times S_1$. Clearly,  $S$ is dominating set. 
We have,  $deg_S((1,0))=|\{0\}\times S_1|=|S_1|=\lceil\frac{|R|}{2}\rceil$  and $deg_{\bar{S}}((1,0))+1=|\bar{S}\cap \{0\}\times R^*|+1=|\{0\}\times R^*|-|S_1|+1=|R|-\lceil\frac{|R|}{2}\rceil\leq \lceil \frac{|R|}{2}\rceil$. Then, $deg_S((1,0))\geq deg_{\bar{S}}((1,0))+1$.
We have, $deg_S((1,x))=|S\cap\{0\}\times Z(R)|=|\{0\}\times Z(R)|-1=|Z(R)^*|$  and $deg_{\bar{S}}((1,x))+1=|\bar{S}\cap \{0\}\times Z(R)|+1=1.$ Hence, $deg_S((1,x))\geq deg_{\bar{S}}((1,x))+1$. 
Let $y\in Z(R)$. If $y^2=0$.  Then, $deg_S((0,y))=|S\cap\Z_2\times \ann_R(y)|=|\{(1,0),(1,x)\}\cup \{0\}\times \ann_R(y)|-2=|\ann_R(y)|$  and $deg_{\bar{S}}((0,y))+1=|\bar{S}\cap \Z_2\times \ann_R(y)|=|\{1\}\times \ann_R(y)|-2+1=|\ann_R(y)|-1$.  Thus, $deg_S((0,y)\geq deg_{\bar{S}}((0,y))$. 
If $y^2\neq 0$,  then $deg_S((0,y))=|S\cap \Z_2\times \ann_R(y)|=|\{(1,0),(1,x)\}|+|\{0\}\times \ann_R(y)|-1=|\ann_R(y)|+1$ 
and $deg_{\bar{S}}((0,y))+1=|\bar{S}\cap \Z_2\times \ann_R(y)|+1=|\{1\}\times \ann_R(y)|-2+1=|\ann_R(y)|-1.$ Then, $deg_S((0,y)\geq deg_{\bar{S}}((0,y))+1$.
Let $u\in R\setminus Z(R)$. Then,  $deg_S((0,u)=1\geq 0+1=deg_{\bar{S}}((0,u))+1$.
Hence, $S$ is a global defensive $1$-alliance.\\
Now, Let $S$ be a global defensive $1$-alliance of minimal cardinality $r= \gamma_1^d(\Gamma(\Z_2\times R))$. 
If $(1,0)\notin S.$ Then, $(0,u)\in S,$ for every $u\in U(R)$ and so $deg_S((0,u))=0\geq deg_{\bar{S}}((0,u))+1=2$  a contradiction. Hence, $(1,0)\in S$ and so $deg_S(1,0)\geq deg_{\bar{S}}((1,0))+1,$ then $|S\cap \{0\}\times R^*|\geq |\bar{S}\cap \{0\}\times R^*|+1.$ Thus,  $|S|\geq \frac{|R|}{2}+1$. 
Suppose that  for all $y\in Z(R)^*$, $(1,y)\in \bar{S}$, then  for all $x\in Z(R)^*$, $(0,x)\in S$.  
 Otherwise, there exists  $(0,y)\in \bar{S}$ which adjacent to at least one element of the form $(0,x)$ ($\Gamma(R)$ is a  connected graph). Thus, $deg_S((0,x))\geq deg_{\bar{S}}(0,x)+1$ and so if $x^2\neq 0$,  we have $|S\cap \Z_2\times \ann_R(x)|\geq |\bar{S} \cap \Z_2 \times \ann_R(x)| +1,$ then  $|S\cap \Z_2\times \ann_R(x)|\geq |\bar{S} \cap \Z_2 \times \ann_R(x)| +1$  and so $|\ann_R(x)|-1\geq |\ann_R(x)|-2$ a contradiction. Thus, for all $x\in Z(R)^*$, $(0,x)\in S$.  Or,  $R$ is a local ring, then there exists $y\in Z(R)^*$ such that $\ann_R(y)=Z(R)$.  Thus, $deg_S((0,y))\geq deg_{\bar{S}}((0,y))+1$ and so  $|S\cap \Z_2\times \ann_R(x)|\geq |\bar{S} \cap \{0\} \times \ann_R(x)| +1$. Then, $|Z(R)^*|\geq |Z(R)^*|+1$   a contradiction.  Hence, there exists $z\in Z(R)^*$ such that $(1,z)\in S$ and so $\lceil\frac{|R|}{2}\rceil+2  \geq |S|\geq \frac{|R|}{2}+2$. 
Thus,  $\gamma_1^d(\Gamma(\Z_2\times R))= \lceil\frac{|R|}{2}\rceil+2$.
\cqfd
Also, we have the following results.

\begin{thm}\label{thmfl}
 Let $R$ be a finite local ring which is not a field. Then,
 \begin{equation*}
 \gamma_{k}^d(\Gamma(\Z_2\times R))=  \left\{
 \begin{array}{ll}
 & 2 \text{\hspace{2.5cm} if } k\in [\![1-|R|;3-|R|]\!],\\&
 \lceil\frac{|R|+k+1}{2}\rceil \text{\hspace{1.1cm} if } k\in [\![4-|R|;4-2|Z(R)|]\!].
 \end{array}
 \right.
 \end{equation*}   
\end{thm}
\pr
\textbf{Case}  $k\in [\![1-|R|;3-|R|]\!]$: Let $x\in Z(R)^*$ such that $Z(R)=\ann_R(x)$ (Since $R$ is a finite local  ring) and set $S=\{(1,0), (0,x)\}$.  Clearly $S$ is a dominating set.  We have,  $deg_S((1,0))=1$ and $deg_{\bar{S}}((1,0))+k=|U(R)|+|Z(R)^*|-1+k=|R|-2+k$.  Then, $deg_S((1,0))\geq deg_{\bar{S}}((1,0))+k$. We have,  $deg_S((0,x))=1$ and $deg_{\bar{S}}((0,x))+k=|Z(R)^*|+|Z(R)^*|-1+k=2|Z(R)|-3+k\leq 2|Z(R)|-|R|$.  So,  $deg_S((0,x))\geq deg_{\bar{S}}((0,x))+k$. Hence, $S$ is a  global defensive $k$-alliance of cardinality $|S|=2$.\\
Now, let $S$ be a global defensive $k$-alliance of minimal cardinality $r=\gamma_{k}^d(\Gamma(\Z_2\times R))$. Then,  $r=|S|\leq 2$ and since $S$ is a dominating set, $|S|\geq 2$. Hence, $r=\gamma_{k}^d(\Gamma(\Z_2\times R))=2$.\\
\textbf{Case}   $k\in [\![4-|R|;4-2|Z(R)|]\!]$: Let $S=\{(1,0),(0,x)\}\cup \{0\}\times S_1$ such that $x\in Z(R)^*$  and  $S_1\subseteq U(R)$ with  $|S_1|=\lceil\frac{|R|+k+1}{2}\rceil -2.$ It is clear that $S$ is a dominating set. We have,  $deg_S((1,0))=1+\lceil\frac{|R|+k+1}{2}\rceil -2= \lceil\frac{|R|+k+1}{2}\rceil-1$ and $deg_{\bar{S}}((1,0))+k= |\{0\}\times R^*|-|S_1|-1+k=|R|-2-\lceil\frac{|R|+k+1}{2}\rceil+k\leq \lceil\frac{|R|+k+1}{2}\rceil-3$. Thus, $deg_S((1,0))\geq deg_{\bar{S}}((1,0))+k$. 
Let $u\in S_1$, then $deg_S((0,u))=1$ and $deg_{\bar{S}}((0,u))+k=k$.  Then, $deg_S((0,u))\geq deg_{\bar{S}}((0,u))+k$. 
We have,  $deg_S((0,x))=1$ and $deg_{\bar{S}}((0,x))+k=2|Z(R)|-3+k\leq 1$. Thus, $deg_S((0,x))\geq deg_{\bar{S}}((0,x))+k$.  Hence, $S$ is a global defensive $k$-alliance of cardinality $|S|=\lceil\frac{|R|+k+1}{2}\rceil$.\\
Now, let $S$ be a global defensive $k$-alliance of minimal cardinality $r=\gamma_{k}^d(\Gamma(\Z_2\times R)).$\\ If $(1,0)\notin S$. Then, $\{0\}\times U(R)\subseteq S$ and also $S$ contains at least one element of the form $(0,z)$ or $(1,z)$ with $z\in Z(R)^*$. Thus, $|S|\geq |U(R)|+1=|R|-|Z(R)|+1$ and so $\frac{|R|+k+1}{2}+1 \leq |S|\leq \lceil\frac{|R|+k+1}{2}\rceil $ a contradiction. Then, $(1,0)$ must be in $S$ and so  $deg_S((1,0))\geq deg_{\bar{S}}((1,0))+k.$ Then, $|S\cup \{0\}\times R^*|\geq |\bar{S}\cup \{0\}\times R^*|+k$ and so $|S\cup \{0\}\times R^*|\geq \frac{|R|+k-1}{2}.$ Then, $|S|\geq \frac{|R|+k+1}{2}$ and by the first part, $|S|\leq \lceil\frac{|R|+k+1}{2}\rceil.$ Hence, $r=|S|=\lceil\frac{|R|+k+1}{2}\rceil$.\\
\cqfd

Notice that for $\Z_4$ when $|Z(\Z_4)|=2$,  $ [\![4-|R|;4-2|Z(R)|]\!]=\{0\}$. In this case, the equalities of Theorems  \ref{thm4.2} and    \ref{thmfl} are the same.

\begin{cor}\label{cor-Z2pn}
	Let $p$ be a prime number and $n$ be positive integer. Then, 
\begin{equation*}	
		\gamma_{k}^d(\Gamma(\Z_2\times \Z_{p^n}))=  \left\{
	\begin{array}{ll}
		& 2 \text{\ \hspace{2cm}  if } k\in [\![1-p^n;3-p^n]\!],\\&
		\lceil\frac{p^n+k+1}{2}\rceil \text{\hspace{0.8cm} if } k\in [\![4-p^n;4-2p^{n-1}]\!], \\&
		\lceil\frac{p^n}{2} \rceil \text{\ \hspace{1.4cm} if } k=-1,\\&
		\lceil \frac{p^n+1}{2} \rceil \text{ \hspace{1cm} if } k=0, \\&
		\lceil\frac{p^n}{2}\rceil+2 \text{\ \hspace{0.75cm} if } k=1.
		
	\end{array}
	\right.
\end{equation*}

and  $\lceil \frac{p^n+k+1}{2} \rceil \leq \gamma_{k}^d(\Gamma(\Z_2\times \Z_{p^n}))\leq \lceil\frac{p^{n-1}(p+2)+k-1}{2}\rceil $   for every $k\in [\![5-2p^{n-1};-2]\!].$	
\end{cor}

For a   local ring $R$   with a nilpotent maximal ideal $M$      of index $2$, we can improve the inequality of Theorem \ref{thm4.1} and give an equality  for the remaining cases other than the ones studied in Theorems  \ref{thm4.2} and   \ref{thmfl}.
We only treat the case where   $|M|\geq 4$. The other cases are rather simples. In fact, if $|M|=2$, then $R\cong \Z_4$ or $\Z_2[X]/(X^2)$  (by \cite[Corollary 1]{MBRB}).  So,    $\gamma_{-3}^d(\Gamma(\Z_2\times R))=\gamma_{-2}^d(\Gamma(\Z_2\times R))=\gamma_{-1}^d(\Gamma(\Z_2\times R))=2$, $\gamma_{0}^d(\Gamma(\Z_2\times R))=3$ and $\gamma_{1}^d(\Gamma(\Z_2\times R))=4$. If $|M|=3$, then $R\cong \Z_9$ or $\Z_3[X]/(X^2)$  (by \cite[Corollary 1]{MBRB}). So,  $\gamma_{-8}^d(\Gamma(\Z_2\times R))=\gamma_{-7}^d(\Gamma(\Z_2\times R))=\gamma_{-6}^d(\Gamma(\Z_2\times R))=2$, $\gamma_{-5}^d(\Gamma(\Z_2\times R))=\gamma_{-4}^d(\Gamma(\Z_2\times R))=3$, $\gamma_{-3}^d(\Gamma(\Z_2\times R))=\gamma_{-2}^d(\Gamma(\Z_2\times R))=4$, $\gamma_{-1}^d(\Gamma(\Z_2\times R))=\gamma_{0}^d(\Gamma(\Z_2\times R))=5$ and $\gamma_{1}^d(\Gamma(\Z_2\times R))=7$.

\begin{thm}\label{thm4.5}
	Let $R$ be a local ring such that its maximal ideal $M$ is nilpotent of index $2$. Then, 
for   $k\in [\![5-2|M|;-2]\!]$ with  $|M|\geq 4$, we have  $\gamma_k^d(\Gamma(\Z_2\times R))=\lceil\frac{|R|+k+1}{2}\rceil$.	
\end{thm}
\pr
Let $k\in [\![5-2|M|;-2]\!]$ and set $S=\{(1,0)\}\cup \{0\}\times S_1\cup \{0\}\times S_2$ with $S_1\subset M^*$ and $S_2\subset U(R)$ such that  $|S_1|=\lceil\frac{2|M|+k+1}{2}\rceil-2$ and $|S_2|=\lceil\frac{|R|+k+1}{2}\rceil-\lceil\frac{2|M|+k+1}{2}\rceil+1$. We have $deg_S((1,0))=|S_1|+|S_2|=\lceil\frac{|R|+k+1}{2}\rceil-1$ and $deg_{\bar{S}}((1,0))+k=|R|-1-(|S_1|+|S_2|)+k=|R|+k+1-\lceil\frac{|R|+k+1}{2}\rceil-1\leq \lceil\frac{|R|+k+1}{2}\rceil-1$. Then, $deg_S((1,0))\geq  deg_{\bar{S}}((1,0))+k$. Let $z\in S_1$, we have  $deg_S((0,z))=|S_1|$ and $deg_{\bar{S}}((0,z))+k=|M|-1-|S_1|+|M|-1+k=2|M|+k-\lceil\frac{2|M|+k+1}{2}\rceil\leq \lceil\frac{2|M|+k+1}{2}\rceil-1$. Then,  $deg_S((0,z))\geq  deg_{\bar{S}}((0,z))+k$. Let $u\in S_2$, we have $deg_S((0,u))=1$ and $deg_{\bar{S}}((0,u))+k=0+k$. Thus,  $deg_S((0,u))\geq deg_{\bar{S}}((0,u))+k$. Hence, $S$ is a global defensive $k-$alliance of cardinality $|S|=\lceil\frac{|R|+k+1}{2}\rceil$. \\
Now, Let $S$ be a global defensive $k-$alliance of minimal cardinality $\gamma_{k}^d((\Z_2\times R))$. From the first part and Theorem \ref{thm4.1} we get the equality $\gamma_{k}^d((\Z_2\times R))=\lceil\frac{|R|+k+1}{2}\rceil$. 
\cqfd   

\begin{cor}
Let $p$ be a prime number. Then, for every $k\in [\![1-p^2;1]\!]$,  $\gamma_k^d(\Gamma(\Z_2\times \Z_{p^2}))=\lceil\frac{p^2+k+1}{2}\rceil.$
\end{cor}

\bigskip

Driss Bennis:  Faculty of Sciences,  Mohammed V University in Rabat,  Morocco.\\
e-mail address:\,driss.bennis@fsr.um5.ac.ma; driss$\_$bennis@hotmail.com.  \\

Brahim El Alaoui: Faculty of Sciences,  Mohammed V University in Rabat,  Morocco.\\
e-mail address:\,brahim.elalaoui@um5r.ac.ma;  brahimelalaoui0019@gmail.com   \\

Khalid Ouarghi:  Department of Mathematics, Faculty of Sciences, King Khalid University. PO Box 960, Abha, Saudi Arabia.\\
e-mail address:\,ouarghi.khalid@hotmail.fr

\end{document}